# ABOUT ADAPTIVE SINGULAR SYSTEMS WITH EXTERNAL DELAY

M. De la Sen

Department of Electricity and Electronics. **Faculty of Science and Technology**

Campus of Leioa (Bizkaia). SPAIN

**Abstract**: This paper is mainly concerned with the robustly stable adaptive control of single-input single-output impulse-free linear time-invariant singular dynamic systems of known order and unknown parameterizations subject to single external point delays. The control law is of pole-placement type and based on input/output measurements and parametrical estimation only. The parametrical estimation incorporates adaptation dead zones to prevent against potential instability caused by disturbances and unmodeled dynamics. The Weierstrass canonical form is investigated in detail to discuss controllability and observability via testable conditions of the given arbitrary state-space realization of the same order.

## I. INTRODUCTION

Singular linear dynamic systems subject to external delays (i.e. in their inputs) have been investigated in [1] by formulating the Weierstrass state-space canonical form. A general formalism for arbitrary state-space realizations is provided in [2]. The so-called Drazin inverse [3] of the singular matrix characterizing the dynamic equation of the system might be used to obtain similarity transformations on the state systems which allow to elucidate if the system is impulsive or not. Instead of using the Drazin inverse and its associate similarity condition, an equivalence transformation may be obtained in the state space through two nonsingular matrices which splits the whole state-space description into a standard dynamic equation plus an algebraic equation involving a nilpotent matrix. Impulsive effects and associate lost of the uniqueness of state- trajectory solution appear if the input is not sufficiently smooth, [1-6]. This paper is devoted to the adaptive control of solvable single-input single-output impulse-free linear time-invariant singular dynamic systems of known orders and unknown parameterizations subject to single external point delays. The system is also potentially subject to disturbances which may potentially include unmodeled dynamics plus external bounded disturbances including noise. The investigation is devoted in several steps. Firstly, the solvability of the dynamic differential- system is investigated for general state-space realizations and also for the Weierstrass canonical form. The trajectory- solution is found under solvability conditions for both the general case and the state-space representation in canonical Weierstrass form. Secondly, the controllability and observability of the Weierstrass canonical form are discussed while testable " ad -hoc" conditions for those properties are obtained. The testability of those conditions is related to the feature that the non-impulsive and impulsive dynamics are described through different equations under the Weierstrass canonical form. The third step consists of synthesizing a pole-placement based controller for the case of known parameters which only requires input/ output measurements. Since Diophantine equations are relevant to the controller parameterization for closed-loop prescribed pole-placement, the system is required to be controllable and observable in order to avoiding zero-pole cancellations in its transfer function. On the other hand, since testing those properties is difficult for arbitrary singular systems, testable controllability/observability conditions are obtained for the Weierstrass canoniocal form of the same order as that of the original state-space realization. This guarantees that the original state-space realization is controllable and observable (and then minimal) under the same conditions so that the transfer function has not zero/pole cancellations. As a result, the Diophantine equation associated with the controller synthesis for closed-loop pole -placement is uniquely solvable and the control law is non- singular if it is sufficiently smooth so that the closed-loop solution is also impulsive-free despite the fact that the transfer function has usually more zeros than poles since the system is singular. The smoothness of the input is achieved if zero-free filters of appropriate minimum orders are used for the input, output and disturbance to be injected in the adaptive control law. The final fourth step is to extend the above controller synthesis to the case of unknown parameterizations by incorporating an adaptive control law which uses parameter estimates instead of true parameters but keeps a similar structure. Also, a relative dead zone is incorporated to prevent against potential instability caused by the presence of disturbances, [7-10]. Such a relative dead zone freezes the parameter estimation at any time if the absolute identification error is sufficiently small compared to a known upper-bounding function of the absolute value of the disturbance. Otherwise, the parameter-adaptive law is kept in operation. The controllability of the estimated model is maintained by using projection of the estimated parameters vector on some convex bounded set where controllability is assumed to be maintained and to which the true parameterization also belongs to. Closed-loop stability is proved under such an adaptive controller.

**I.1. Notation**: **N, R** and **C** are, respectively, the sets of natural, real and complex numbers, $\mathbf{N}_0 := \mathbf{N} \cup \{0\}$, $\mathbf{R}_0^+ := \{z \in \mathbf{R}: z \geq 0\}$ and $\mathbf{C}_0^+ := \{z \in \mathbf{C}: \operatorname{Re} z \geq 0\}$

$\operatorname{Lap}(.)$ and $\operatorname{Lap}^{-1}(.)$ denote, respectively, the unilateral left Laplace transform and corresponding anti-transform provided that they exist such that $V(s) = \operatorname{Lap}(v(t))$, $v(t) = \operatorname{Lap}^{-1}(V(s))$ are corresponding transform / anti-transform pairs (roughly speaking, low-case and its corresponding capital stand for a certain signal and its Laplace transform). Since the Laplace transform argument "s" is (formally) equivalent to the time- derivative operator $D := \frac{d}{dt}$ (extended to $D^i := \frac{d^i}{dt^i}$ with $D^0 \equiv 1$ and $D^1 \equiv D$) then $V(s) = \operatorname{Lap}(v(t))$, if it exists, is formally equivalent to $V(D)$ provided that v(t) is an original function; i.e. $v(t) = 0$; $\forall t \leq 0$.

I denotes the identity matrix which may be endowed with a subscript defining its order if such an order is relevant.

$C^{(q)}(D_R, \mathbf{R}^p)$ is the space of vector real functions of class q of definition domain $D_R \subset \mathbf{R}$ for any integer $p \geq 1$. If $v \in C^{(q)}(D_R, \mathbf{R})$ then $v^{(k)}(t) = D^{(k)}v(t)$ for any integer k subject to $0 \leq k \leq q$.

## II. THE SINGULAR CONTROLLED PLANT

Consider the SISO linear time-invariant dynamical system

$$E\dot{x}(t) = Ax(t) + bu(t) + du(t-h) + v_0(t)$$
$$y(t) = c^T x(t)$$
$$x(0) = x_0, \; u(t) = \Psi(t), \; \forall t \in [-h, 0] \quad (1)$$

where $x(t) \in \mathbf{R}^n$ is the state-vector, $u(t) \in \mathbf{R}$, $y(t) \in \mathbf{R}$ and $v_0(t) \in \mathbf{R}$ are the control input (being differentiable $\ell - 1$ times with $\ell$ being the index of E, [2-3]), output and disturbance, respectively; $E \in \mathbf{R}^{n \times n}$ is a singular matrix, $A \in \mathbf{R}^{n \times n}$, $b \in \mathbf{R}^n$, $d \in \mathbf{R}^n$, $c \in \mathbf{R}^n$; $h > 0$ is a control delay and $\Psi(t) \in \mathbf{R}^n$ is the initial control function. The disturbance may include contributions of unmodeled dynamics and bounded external disturbances. If the disturbance $v_0 \equiv 0$ on $\mathbf{R}_0^+$ then the system (1) will be said to be the nominal system. The following assumption is made through the paper, [1-2]:

**Assumption 1**. The pair $(E, A)$ is solvable; i.e. the matrix pencil $(E + sA)$ is regular, namely, $\text{Det}(E - sA)$ is not identically zero for $\forall s \in \mathbf{C}$. □

Some especial particular results through the paper will be obtained when, in addition, the subsequent assumption holds as well, [4]:

**Assumption 2**. The pair $(E, A)$ is impulse-free, namely, $\text{rank}(E) = \deg(\text{Det}(E - sA))$. □

If Assumption 1 holds there exist, in general, non-unique state trajectory solutions. Those solutions become unique for all time within some appropriate manifolds containing the so-called admissible initial state sets, which depend on initial values for the input vector, provided that the bounded initial conditions belong to such a manifold and that such an input vector is sufficiently smooth. If the input is not sufficiently smooth then the solution is impulsive and then non-unique even if Assumption 1 holds. Assumption 1 is necessary to have a well-posed problem. If Assumptions 1-2 hold together, then a unique (and then non-impulsive) solution exists for all time for each bounded initial conditions constrained within the above mentioned manifold and each everywhere piecewise continuous control vector function, [2], [4]. It has been proved in [2] that Assumption 1 holds if and only if $Ker(E \cap X_i) = 0$ for i =0, 1, 2, ... , where $X_0 := Ker(A)$ ; $X_i := \{x : Ax \in EX_{i-1}\}$, $\forall i \in \mathbf{N}$, and, equivalently, if and only if

$$\text{rank} \underbrace{\begin{bmatrix} E & 0 & \cdots & 0 \\ A & E & 0 \cdots & 0 \\ 0 & A & E \cdots & 0 \\ & & \ddots & E \\ & & & A \end{bmatrix}}_{n} \Bigg\} n+1 = n^2$$

If $\ell := ind(E) \in \mathbf{N}_0 := \mathbf{N} \cup \{0\}$, the index of E, then there is a nilpotent matrix $N_0$ with $\ell := n\,ind(N_0)$, the nilpotency index of E, such that, [3],

$$\ell := ind(E) = Min\left(k \in \mathbf{N}_0 : rank(E^k) = rank(E^{k+1})\right) = Min\left(k \in \mathbf{N}_0 : Ker(E^k) = Ker(E^{k+1})\right)$$
$$= n\,index(N_0) := Min\left(k \in \mathbf{N}_0 : E^k = 0\right)$$

$$E = T^{-1} Diag(E_0, N_0) T \; ; \; E^D = T^{-1} Diag(E_0^{-1}, 0) T$$

where T is a nonsingular real n-matrix, $E^D$ is the Drazin inverse of E and $E_0 \in \mathbf{R}^{n_1 \times n_1}$ is nonsingular with $n_1 < n$. If Assumption 2 holds then trivially $N_0 = 0$; i.e. the idempotent matrix becomes zero. The following auxiliary result extends directly a previous one proved in [2], and quoted from [1] for the disturbance-free case, by considering the singular system (1) with incorporates additive disturbances in the state dynamics as well as an output equation. It is concerned with the equivalence of the state-space description (1) with an associate Weierstrass canonical form:

**Lemma 1**. If Assumption 1 holds (i.e. if the pair $(E, A)$ is solvable) then there exist nonsingular real matrices P and Q of order $n \times n$ such that the state-space description (1) is equivalent to the Weierstrass canonical form:

$$\dot{z}_1(t) = W z_1(t) + \alpha_1 u(t) + \beta_1 u(t-h) + \eta_1(t) \tag{2.a}$$
$$N \dot{z}_2(t) = z_2(t) + \alpha_2 u(t) + \beta_2 u(t-h) + \eta_2(t) \tag{2.b}$$
$$y(t) = c^T P z(t) = \gamma^T z(t) = \gamma_1^T z_1(t) + \gamma_2^T z_2(t) \tag{2.c}$$

where

$$z(t) = (z_1^T(t), z_2^T(t))^T = P^{-1} x(t) \in \mathbf{R}^n \text{ with } z_i \in \mathbf{R}^{n_i} \ (i=1,2)$$

$$QEP = \text{Diag}(I_{n_1}, N) \ ; \ QAP = \text{Diag}(W, I_{n_2})$$

$$Qb = (\alpha_1^T, \alpha_2^T)^T; \ Qd = (\beta_1^T, \beta_2^T)^T; \ Qv_0 = (\eta_1^T, \eta_2^T)^T \ ; \ \gamma := Pc = (\gamma_1^T, \gamma_2^T)^T$$

(2.d)

with $\alpha_i, \beta_i, \gamma_i, \eta_i \in \mathbf{R}^{n_i} \ (i=1,2)$ and W being a nonsingular real matrix of order $n_1 < n$ and N is a nilpotent real matrix of order $n_2 = n - n_1$ with $n\text{index}(N) = \ell = \text{ind}(E)$. Assume that $u \in C^{(\ell-1)}(\mathbf{R}_0^+ \cup [-h, 0], \mathbf{R})$ and $\eta_2 \in C^{(\ell-1)}(\mathbf{R}_0^+ \cup [-h, 0], \mathbf{R}^{n_2})$ with $u(t) = \Psi(t)$ and $\eta_2 \in C^{(\ell-1)}(\mathbf{R}_0^+ \cup [-h, 0], \mathbf{R}^{n_2})$, $\forall t \in [-h, 0]$ for some $\Psi \in C^{(\ell-1)}([-h, 0], \mathbf{R})$ and $\Psi_{v02}(t) \in C^{(\ell-1)}([-h, 0], \mathbf{R}^{n_2})$. Then, the unique state-trajectory solution of the Weierstrass canonical form $\mathbf{R}_0^+$ is

$$z_1(t) = e^{Wt} \left[ z_{10} + \int_0^t e^{-W\tau} (\alpha_1 u(\tau) + \alpha_2 u(t-\tau) + \eta_1(\tau)) d\tau \right] \quad (3.a)$$

$$z_2(t) = -\sum_{i=0}^{\ell-1} N^i \left[ \alpha_2 u^{(i)}(t) + \beta_2 u^{(i)}(t-h) + \eta_2(t) \right] \quad (3.b)$$

for $t \geq h$ with $z_{10} = z_1(0)$, and

$$z_1(t) = e^{Wt} \left[ z_{10} + \int_0^t e^{-W\tau} (\alpha_1 u(\tau) + \alpha_2 \Psi(t-\tau) + \eta_1(\tau)) d\tau \right] \quad (4.a)$$

$$z_2(t) = -\sum_{i=0}^{\ell-1} N^i \left[ \alpha_2 \Psi^{(i)}(t) + \beta_2 \Psi^{(i)}(t-h) + \eta_2^{(i)}(t) \right] \quad (4.b)$$

for $t \in [0, h)$. The unique output trajectory solution on $\mathbf{R}_0^+$ is given by combining (2.c) with (3) and (4). If $\ell = 0$ (the system is standard, i.e. nonsingular) or $\ell = 1$ then the state and output trajectories are always unique on $\mathbf{R}_0^+$. Assume that $\ell \geq 2$ and $u \in C^{(j)}(\mathbf{R}_0^+ \cup [-h, 0], \mathbf{R})$ for some $0 \leq j \leq \ell - 2$ but $u \notin C^{(\ell-1)}(\mathbf{R}_0^+ \cup [-h, 0], \mathbf{R})$ then the state and output trajectory solutions on $\mathbf{R}_0^+$ are non-unique and $u^{(j+k)}(t)$ is replaced with $\delta^{(j+k)}(t)$ (the k-th order Dirac distribution) everywhere where the (high-order) time derivative does not exist.

If Assumptions 1-2 hold (i.e. if the pair $(E, A)$ is solvable and impulse-free) then $N \equiv 0$ so that the second equation is purely algebraic for all time taking the form:

$$z_2(t) = -(\alpha_2 u(t) + \beta_2 u(t-h) + \eta_2(t)) \quad ; \quad \forall t \geq h$$

$$z_2(t) = -(\alpha_2 u(t) + \beta_2 \Psi(t-h) + \eta_2(t)) \quad ; \quad \forall t \in [0, h)$$

subject to $u(t) = \Psi(t), \forall t \in [-h, 0]$.

**Outline of Proof**: The existence of the state-space description (2) under the given conditions is proved in [1-2]. The remaining part of the proof including the form of the solutions and their potential uniqueness follows directly from (2) and the use of high-order distributional derivatives for singular impulsive time-delay systems, [5]. □

The following direct general result concerning the description of (1) by a time-differential input-output model under Assumption 1 is proved in Appendix A:

**Theorem 1**. The following properties hold:

**(i)** If Assumption 1 holds then the singular system (1) is described by the differential equation:

$$M(D)y(t) = \Delta_0(D)u(t) + \Delta_1(D)u(t-h) + \Delta_2^T(D)v_0(t) \tag{5}$$

provided that $D^i u(t)$ $(i = 0, 1, ..., \ell - 1)$ exists subject to initial conditions $D^i y(0)$ $(i = 0, 1, ..., \ell - 1)$, where

$$M(D) := \text{Det}(DE - A) \quad ; \quad \Delta_0(D) := \Delta_2^T(D)b \tag{6.a}$$

$$\Delta_1(D) := \Delta_2^T(D)d \quad ; \quad \Delta_2^T(D) := c^T \text{Adj}(DE - A) \tag{6.b}$$

Assume that identical (zero-free and zero/pole cancellation free) linear filters order of transfer functions $G_f(s) = 1/F(s)$, with $F(s)$ being a monic polynomial of arbitrary degree $n_F \geq 1$, are used to filter the various signals of the system (1) so that the filtered equations

$$U_f(s) = \frac{1}{F(s)}U(s), \; Y_f(s) = \frac{1}{F(s)}Y(s), \; V_{0f}(s) = \frac{1}{F(s)}V_0(s), \; X_{0f}(s) = \frac{1}{F(s)}(Ex_0) \tag{7}$$

have respective bounded initial conditions $x_{uF0} = x_{uF0}(0)$, $x_{yF0} = x_{yF0}(0)$, $x_{0F0} = x_{0F}(0)$ and $x_{v0F0} = x_{v_0 F0}(0)$ in $\mathbf{R}^n$. If Assumption 1 holds then there exist solutions of the differential equation (2) that satisfy

$$y(t) = (F(D) - M(D))y_f(t) + \Delta_0(D)u_f(t) + \Delta_1(D)u_f(t-h) + \Delta_2^T(D)v_{0f}(t) + y_0(t)$$

(8)

for any set of bounded initial conditions of the filters, where

$$y_0(t) := M(D)\left\{G'^T(D)[v_{0f}(t) - G_f(D)x_{v0F0}] - G_f'^T(D)G(D)x_{uF0}\right\}$$
$$+ G_f'^T(D)(M(D)-1)x_{yF0} \qquad (9)$$

where

$$G(s) := c^T(sE - A)^{-1}(b + de^{-hs}) = M^{-1}(s)(\Delta_0(s) + \Delta_1(s)e^{-hs}) \qquad (10)$$

is the transfer function of the system (1) related to the control input and $y_0(t)$ is a signal being dependent of the external disturbance and the initial conditions of the various filters. □

If the closed-loop system is stabilizable; i.e. $\text{rank}[sE - A, b + de^{-hs}] = n$, $\forall s \in \mathbf{C}_0^+$, and it is stabilized via some linear controller and its associate stabilizing control law then Theorem 1 applies to the closed-loop system with the appropriate replacements for the closed-loop parameterization and feedback control input, namely, $M(s) \to M_m(s)$, $\Delta_i(s) \to \Delta_{mi}(s)$ (i=1, 2) and $u(t) \to r(t)$ (the external reference input). If, in addition, the polynomial F(s) is Hurwitz, i.e. the various filters are asymptotically stable, what is standard in engineering applications; then $y_0(t) \to 0$ as $t \to \infty$ for any bounded initial conditions of the filters. As a result, such a signal may be neglected in both the non-adaptive and adaptive controller synthesis process for stabilization purposes. It may also be removed from the estimation equations in the adaptive case as it will be discussed in Section IV. The following result is devoted to the uniqueness of the state and output trajectories of (1) for all time if the input is sufficiently smooth and the initial conditions belong to some appropriate set, even if Assumption 2 does not hold, or if Assumption 2 holds and the initial state is in some appropriate subset of the whole state space irrespective of the input being sufficiently smooth or not. Such a result corresponds to a parallel one obtained in [5] for the case of singular time-delay systems with multiple internal (i.e. in the state) single constant point delays which investigates the whole set of potentially impulsive solutions for the case where the input is not necessarily sufficiently smooth and Assumption 2 does not necessarily hold.

**Theorem 2.** Assume that $u \in C^{(\ell-1)}(\mathbf{R}_0^+ \cup [-h, 0), \mathbf{R})$ and the external disturbance of the form $v_0(t) = Q^{-1}(\eta_1^T(t), \eta_2^T(t))^T$ satisfies $\eta_2 \in C^{(\ell-1)}(\mathbf{R}_0^+ \cup [-h, 0), \mathbf{R}^{n_2})$ with $u(t) = \Psi(t)$ and $\eta_2(t) = \Psi_{v02}(t)$ for all $t \in [-h, 0)$ for some given $\Psi \in C^{(\ell-1)}([-h, 0], \mathbf{R}^{n_2})$ and $\Psi_{v02} \in C^{(\ell-1)}([-h, 0], \mathbf{R}^{n_2})$ and, furthermore,

Assumption 1 holds. Define the sets of admissible initial states $S(\Psi) \subset \mathbf{R}^n$ and admissible initial data $A(\Psi) \subset \mathbf{R}^n \times C^{(\ell-1)}([-h,0], \mathbf{R}) \times C^{(\ell-1)}([-h,0], \mathbf{R}^{n_2})$, respectively, as follows:

$$S(\Psi) := \left\{ z = (z_1^T, z_2^T)^T, z_1 \in \mathbf{R}^{n_1}, z_2 = -\sum_{i=0}^{\ell-1} N^i [\alpha_2 u^{(i)}(0) + \beta_2 \Psi^{(i)}(-h) + \Psi_{v02}^{(i)}(0)] \in \mathbf{R}^{n_2} \right\}$$

$$A(\Psi) := \left\{ (z_0, \Psi, \Psi_{v02}) : z_0 \in S(\Psi), \Psi \in C^{(\ell-1)}([-h,0], \mathbf{R}), \Psi_{v02} \in C^{(\ell-1)}([-h,0], \mathbf{R}^{n_2}) \right\}$$

If $\ell = 0,1$ or if $(z_0, \Psi:[-h,0] \to \mathbf{R}, \Psi_{v02}:[-h,0] \to \mathbf{R}^{n_2}) \in A(\Psi)$ with $\ell \geq 2$ then the unique output trajectory solution to (1) on $\mathbf{R}_0^+$ is defined by

$$y(t) = \gamma_{11} e^{Wt} \left[ z_{10} + \int_0^t e^{-W\tau} (\alpha_1 u(\tau) + \alpha_2 u(t-\tau) + \eta_1(\tau)) d\tau \right]$$
$$- \gamma_{12} \sum_{i=0}^{\ell-1} N^i [\alpha_2 u^{(i)}(t) + \beta_2 u^{(i)}(t-h) + \eta_2^{(i)}(t)] ; \forall t > h \qquad (11.a)$$

$z_{10} = z_1(0)$ subject to

$$y(t) = \gamma_{11} e^{Wt} \left[ z_{10} + \int_0^t e^{-W\tau} (\alpha_1 u(\tau) + \alpha_2 \Psi(t-\tau) + \eta_1(\tau)) d\tau \right]$$
$$- \gamma_{12} \sum_{i=0}^{\ell-1} N^i [\alpha_2 u^{(i)}(t) + \beta_2 \Psi^{(i)}(t-h) + \Psi_{v02}^{(i)}(t)(t)] ; \forall t \in [0,h] \qquad (11.b)$$

where

$$\gamma_i = c_1^T P_{1i} + c_2^T P_{2i} ; \quad \alpha_i = Q_{1i} b_1 + Q_{2i} b_2 \qquad (12.a)$$

$$\beta_i = Q_{1i} d_1 + Q_{2i} d_2 ; \quad \eta_i = Q_{1i} v_{01} + Q_{2i} v_{02} \qquad (12.b)$$

for $i = 1,2$ where $P = \text{Block Matrix}(P_{ij}; i,j = 1,2)$, $Q = \text{Block Matrix}(Q_{ij}; i,j = 1,2)$, $\gamma = (\gamma_1^T, \gamma_2^T)^T$; $\alpha = (\alpha_1^T, \alpha_2^T)^T$; $\beta = (\beta_1^T, \beta_2^T)^T$; $v_0 = (v_{01}^T, v_{02}^T)^T$ with $P_{ij}, Q_{ij} \in \mathbf{R}^{n_i \times n_i}$, ; $\gamma, \alpha, \beta, v_0 \in \mathbf{R}^{n_i}$. In particular, the output trajectory solution becomes

$$y(t) = \gamma_{11} e^{Wt} \left[ z_{10} + \int_0^t e^{-W\tau} (\alpha_1 u(\tau) + \alpha_2 u(t-\tau) + \eta_1(\tau)) d\tau \right]$$
$$- \gamma_{12} [\alpha_2 u(t) + \beta_2 u(t-h) + \eta_2(t)] ; \forall t > h \qquad (13.a)$$

$$y(t) = \gamma_{11} e^{Wt} \left[ z_{10} + \int_0^t e^{-W\tau} (\alpha_1 u(\tau) + \alpha_2 \Psi(t-\tau) + \eta_1(\tau)) d\tau \right]$$
$$- \gamma_{12} [\alpha_2 u(t) + \beta_2 \Psi(t-h) + \eta_2(t)] ; \forall t \in [0,h] \qquad (13.b)$$

for $\ell = 1$ (singular solvable system with zero nilpotent matrix) and $\left( z_0, \Psi:[-h,0] \to \mathbf{R}, \Psi_{v02}:[-h,0] \to \mathbf{R}^{n_2} \right) \in A(\Psi)$ even if $u \notin C^{(\ell-1)}\left( \mathbf{R}_0^+ \cup [-h,0), \mathbf{R} \right)$ or $\Psi_{v02} \notin C^{(\ell-1)}\left( \mathbf{R}_0^+ \cup [-h,0), \mathbf{R}^{n_2} \right)$; and

$$y(t) = \gamma_{11} e^{Wt} \left[ z_{10} + \int_0^t e^{-W\tau} \left( \alpha_1 u(\tau) + \alpha_2 u(t-\tau) + \eta_1(\tau) \right) d\tau \right] \quad ; \forall t > h$$

(14.a)

$$y(t) = \gamma_{11} e^{Wt} \left[ z_{10} + \int_0^t e^{-W\tau} \left( \alpha_1 u(\tau) + \alpha_2 \Psi(t-\tau) + \eta_1(\tau) \right) d\tau \right] \quad ; \forall t \in [0,h]$$

(14.b)

for $\ell = 0$ (nonsingular system ; i.e. system in standard form such that Assumption 1 always holds) and $\left( z_0, \Psi:[-h,0] \to \mathbf{R}, \Psi_{v02}:[-h,0] \to \mathbf{R}^{n_2} \right) \in A(\Psi, \Psi_{v02})$ even if $u \notin C^{(\ell-1)}\left( \mathbf{R}_0^+ \cup [-h,0), \mathbf{R} \right)$ or $\Psi_{v02} \notin C^{(\ell-1)}\left( \mathbf{R}_0^+ \cup [-h,0), \mathbf{R}^{n_2} \right)$. □

The proof of Theorem 2 follows directly by combining Lemma 1 and (2.c) since any existing output trajectory solution is independent of the state-space representation Note that for the cases $\ell = 0$ (nonsingular system also referred to as system in standard form, [6]) and $\ell = 1$ (singular system with identically zero idempotent matrix) there exists a unique output trajectory solution for each set of admissible data irrespective of the input being sufficiently smooth or not. In addition, a system in standard system always satisfies Assumption 1; i.e. the pair $(I_n, A)$ is always solvable. See [5-6] for the solutions of the disturbance- free case; i.e. $v_0 \equiv 0$, under multiple internal point delays. Note also that the external disturbance relation in-between the Weiertrass canonical form to the state-space representation (1) by $v_0(t) = Q^{-1} \left( \eta_1^T(t), \eta_2^T(t) \right)^T$ is uniquely defined by an associate bijective mapping since Q is nonsingular. Therefore, the constraint $\eta_2(t) = \Psi_{v02}(t)$ for all $t \in [-h,0]$ required by the set of admissible data is equivalent to $v_0(t) = \left( \Psi_1^T(t), \Psi_1^T(t) \right)^T$ with any arbitrary $\Psi_2 \in C^{(\ell-1)}\left( [-h,0], \mathbf{R}^{n_2} \right)$ for all $t \in [-h,0]$ and $\Psi_1[-h,0] \to \mathbf{R}^{n_1}$ being absolutely continuous with isolated bounded discontinuities on a set of zero measure in order to guarantee the uniqueness of the output trajectory solution on $\mathbf{R}_0^+$ from the well-known Picard-Lindelöf theorem and the properties of the set of admissible data for the solvable singular system.

**Remark 1**. The controllability and observability of the Weiertrass canonical form are investigated in Appendix B where Theorem B.1 provided with specific related tests. The controllability and observability properties will be relevant through the subsequent sections of the manuscript for the controller synthesis problem which involves Diophantine Equations of polynomials On the other hand, an alternative expression for the output to those given in Eqs. (5)- (6) and in Eqs. (A.1), subject to (7), is given in Eqs. B.6 - B.7 in Appendix B, without using filtered signals, based on the Weierstrass canonical form by taking advantage of the independence of the output trajectory solution of each particular state-space representation. □

### III. CONTROL LAW AND STABILIZING CONTROLLER

**A. Case of known parameters**

Throughout this section the time –derivative and Laplace operators D and s are used indistinctly depending on the physical context, namely, D is used in differential equations and s in transfer functions or matrices and associated numerator and denominator polynomials or quasi-polynomials. The control law is generated by a stabilizing controller which is built with a feedforward compensator (precompensator), a feedback compensator of respective transfer functions $H_{ff}(s) = \frac{T(s)}{R(s)}$ and $H_{fb}(s) = \frac{S(s)}{R(s)}$, the filters of transfer function $1/F(s)$ and an smooth disturbance compensating signal $L(D)u_0(t)$ and it is defined implicitly by

$$R(D,h)u(t) = T(D)u_{cf0}(t) - S(D,h)y_{f0}(t) + L(D)u_{0f0}(t) \qquad (15)$$

where $u_c(t)$ is a bounded reference signal the various filtered signals (subscripted with "fo") are obtained through the filters of transfer functions $1/F_0(s)$ using the corresponding unfiltered signals as inputs. The filters of transfer functions $1/F_0(s)$ might optionally be distinct of those of transfer functions $1/F_0(s)$ referred to in Theorem 1 for generalization purposes. $T(s)$ is a polynomial and $R(s,h)$ and $S(s,h)$ are, in general, dependent on the delay quasi-polynomials:

$$R(s,h) = R_0(s) + e^{-hs}R_1(s) \; ; \; S(s,h) = S_0(s) + e^{-hs}S_1(s) \qquad (16)$$

and $R_i(s)$ and $S_i(s)$ (i =0,1) are polynomials. The degrees of the involved polynomials are denoted as:

$$n_M = \deg(M(s)), \; n_{\Delta_i} = \deg(\Delta_i(s)), \; n_L = \deg(L(s)) \qquad (17.a)$$

$$n_{F_0} = \deg(F_0(s)), \; n_T = \deg(T(s)), \; n_{R_i} = \deg(R_i(s)), \; n_{S_i} = \deg(S_i(s)) \qquad (17.b)$$

for i=0,1. The main technical reason to involve the use of filtered signals is the need to measure time-derivatives of the various involved functions in both the non-adaptive and later adaptive versions. Equivalently, (15) may be expressed with the unfiltered signals as follows:

$$F_0(D)R(D,h)u(t) = T(D)u_c(t) - S(D,h)y(t) + L(D)u_0(t) \tag{18}$$

The closed-loop system time-differential equation, obtained by combining (5) and (18), is

$$\begin{aligned}
&(F_0(D)M(D)R_0(D) + \Delta_0(D)S_0(D))y(t) \\
&+ (F_0(D)M(D)R_1(D) + \Delta_0(D)S_1(D) + \Delta_1(D)S_0(D))y(t-h) + \Delta_1(D)S_1(D)y(t-2h) \\
&= M_0^*(D)y(t) + M_1^*(D)y(t-h) + T(D)(\Delta_0(D)u_c(t) + \Delta_1(D)u_c(t-h)) + \sigma_0(t)
\end{aligned} \tag{19}$$

where

$$M_0^*(s) = F_0(s)M(s)R_0(s) + \Delta_0(s)S_0(s) \tag{20.a}$$

$$M_1^*(s) - \Delta_1(s)S_0(s) = F_0(s)M(s)R_1(s) + \Delta_0(s)S_1(s) \tag{20.b}$$

$$\begin{aligned}
\sigma_0(t) &= L(D)(\Delta_0(D)u_0(t) + \Delta_1(D)u_0(t-h)) \\
&+ F(D)M(D)(R_0(D)\Delta_2^T v_0(t) + R_1(D)\Delta_2^T v_0(t-h)) - \Delta_1(D)S_1(D)y(t-2h)
\end{aligned} \tag{20.c}$$

with

$$n_{M_0^*} = \deg(M_0^*(s)) \ , \ n_{M_1^*} = \deg(M_1^*(s)) \tag{21}$$

The subsequent result related to closed-loop stability and model matching property is proved in Appendix A.

**Theorem 3**. The following properties hold:

**(i)** Assume that the system (1) is solvable and controllable and observable and that the external disturbance $v_0 : C^{(\ell-1)}(\mathbf{R}_0^+, \mathbf{R}^n)$ is uniformly bounded. Then, the closed-loop system (19) - (20), subject to (21), is globally Lyapunov´s stable if the following conditions hold:

1. The polynomial F(s) is Hurwitz and has no common zero with $\Delta(s)$. Also, the quasi-polynomial $M^*(s) = M_0^*(s) + e^{-hs}M_1^*(s)$ is built with prescribed monic polynomials $M_i^*(s)$ subject to degree constraints $n_{M_i^*} = 2(n_F + n_M) - 1$ (i=0,1) while being Hurwitz with all its zeros in $\mathbf{C}_{-\nu} := \{z \in \mathbf{C} \setminus \mathbf{C}_0^+ : \operatorname{Re} z \leq -\nu\}$ for some $\nu \in \mathbf{R}_0^+$.

2. The control law (15)–(16) is generated with polynomials $R_i(s)$ (monic) and $S_i(s)$ ($i=0,1$) being respective solutions to (20.a)- (20.b) which exist and are unique subject to degrees

$$n_{R_i} = n_{F_0} + n_M - 1; \quad n_{S_i} = n_{F_0} + n_M - 2 \text{ for } i=0,1 \tag{22.a}$$

provided that

$$n_L \geq \text{Max}\left(0, n_{F_0} + n_M - 2\right), \quad n_{F_0} \geq \text{Max}\left(0, \frac{n_T + \ell - n_M}{2}\right) \tag{22.b}$$

3. $\displaystyle \sup_{\text{Re } s = v_1} \left( \left| \frac{\Delta_1(s) S_1(s)}{M_0^*(s) e^{2hs} + M_1^*(s) e^{hs}} \right| \right) < 1, \quad \forall s \in \mathbf{C} \text{ fulfilling } \text{Re } s = v_1 \in (0, v) \text{ for some } v_1 \in \mathbf{R}^+$.

**(ii)** Assume, in addition, that the quasi-polynomial $\Delta(s) = \Delta_0(s) + \Delta_1(D) e^{-hs}$ is Hurwitz and that the compensating signal $u_0(t)$ is generated as:

$$u_0(t) = \int_0^t h_{0v}^T(\tau) v_0(t-\tau) d\tau + h_{0y}(\tau) y(\tau) d\tau \tag{23}$$

with impulse responses $h_{0v}: \mathbf{R}_0^+ \to \mathbf{R}^n$ and $h_{0y}: \mathbf{R}_0^+ \to \mathbf{R}$ associated with respective stable transfer matrix and function:

$$H_{0v}(s) = \frac{1}{L(s)\Delta(s)} \left( F_0(s) M(s) R(s) \Delta_2^T(s) \left(1 + e^{-hs}\right) \right); \quad H_{0y}(s) = \frac{\Delta_1(s) S_1(s)}{L(s)\Delta(s)} e^{-2hs}$$

$$\tag{24}$$

Then, asymptotic model–matching of any reference model of transfer function $G_m(s) = \frac{\Delta(s) T(s)}{M_m(s)}$ with $M_m(s) = M_0^*(s) + M_1^*(s) e^{-hs} + \Delta_1(s) S_1(s) e^{-2hs}$ being Hurwitz and satisfying the degree realizability constraint $n_m \geq n - 1 + n_T$ is achieved for a reference signal $u_c: \mathbf{R}_0^+ \to \mathbf{R}$ under any given set of bounded initial conditions.

**(iii)** Assume that the conditions of Property (i) hold, $\Delta_1(s) \equiv 0$ and the external disturbance $v_0 \equiv 0$ on $\mathbf{R}_0^+$. Assume also that the factorization $\Delta(s) = \Delta_0(s) = \Delta_0^+(s) \Delta_0^-(s)$ is performed with $\Delta_0^+(s)$ being monic with zeros, if any (otherwise $\Delta_0^+(s) = 1$), within $\mathbf{C}_{-v}$ and $\Delta_0^-(s)$ having all its zeros within $\mathbf{C} \setminus \mathbf{C}_{-v}$, for some prescribed $v \in \mathbf{R}_+$. Then, the stable reference model of transfer function $G_m(s) = \dfrac{\Delta_m(s)}{M_m(s)} = \dfrac{\Delta_0^-(s)\Delta_m'(s)}{M_m(s)}$ is matched for all time under zero initial conditions of both the reference model and the system (1) by a control law (18)-(20) with the following restrictions:

$$R_1(s) = S_1(s) = M_1^*(s) = 0\;;\; T(s) = \Delta_m'(s)\;;\; R_0(s) = \Delta_0^+(s) R_0'(s)\;;\; u_0 \equiv 0 \text{ on } \mathbf{R}_0^+$$

(25.a)

$$M_m(s) = M_0^*(s) = F_0(s) M(s) R_0'(s) + \Delta_0^-(s) S_0(s) \quad (25.b)$$

$$n_{R_0'} = n_{F_0} + n_M - n_{\Delta_0} - 1\;;\; n_{S_i} = n_{F_0} + n_M - n_{\Delta_0} - 2 \quad (25.c)$$

with $M_m(s) = M_0^*(s)$ and $F_0(s)$ $R_0'(s)$ being monic. □

The fact that $M_i^*(s)$ and $R_i(s)$ (i=0,1) are assumed monic is only relevant for the uniqueness of the solution to the Diophantine equations (19.a)-(19.b) but not for the stability and model matching properties. The transfer function of the disturbance–free (nominal) system related to the control input is zero-pole cancellation-free if and only if the realization of order n is controllable and observable which may be directly tested by using the Weierstrass canonical state space description via Theorem 2. The following result is a particular case of Theorem 3 for the case when the feed-forward and feedback compensator transfer functions are $H_{ff}(s) = \dfrac{T(s)}{R_0(s)}$ and $H_{fb}(s) = \dfrac{S_0(s)}{R_0(s)}$, respectively. Its proof follows directly under similar arguments as those invoked to prove Theorem 3 in Appendix A.

**Corollary 1**. Assume that the system (1) is solvable, controllable and observable, the external disturbance $v_0 : \mathbf{R}_0^+ \to \mathbf{R}^n$ is uniformly bounded, $M_0^*(s)$ is Hurwitz with its zeros in $\mathbf{C}_{-\nu}$ subject to the degree constrain of Theorem 3, $M_1^*(s) \equiv 0$ and $F_0(s)$ has no common zero with $\Delta(s)$. Assume also that the control law generated from (15)-(16) with polynomials $R_1(s) = S_1(s) \equiv 0$ and polynomials $R_0(s)$ and $S_0(s)$ which exist and are unique solutions to (20.a) under the degree constraint of Theorem 3. Then, the particular closed-loop system resulting from (19)-(20) is globally Lyapunov´s stable if $\underset{\operatorname{Re} s = \nu_1}{\operatorname{Sup}} \left( \left| \dfrac{\Delta_0(s) S_1(s)}{M_0^*(s) e^{hs}} \right| \right) < 1,\; \forall s \in \mathbf{C}$ fulfilling $\operatorname{Re} s = \nu_1 \in (0, \nu)$ for some $\varepsilon_1 \in \mathbf{R}^+$. □

**Remark 2**. Note that the index of E which equalizes the nilpotency index $\ell$ of the nilpotent matrix arising in the Weierstrass decomposition has to be taken into account in the design of the control law by using a filter of appropriate order of transfer function $1/F_0(s)$ since the control has to be $(\ell-1)$- times time-differentiable for uniqueness of the solution of (1) under the appropriate initial conditions. However, since $\ell \leq n_2 \leq n$, it would be necessary to increase the degree of $F_0(s)$ in the adaptive case in an amount $\text{Max}(n-n_M, n-\ell+1)$ satisfying:

$$\text{Max}(n-n_M, n-n_2) \leq \text{Max}(n-n_M, n-\ell+1) \leq \text{Max}(n-n_1, n-\ell+1)$$

in order to guarantee that $y^{(j)}(t)$, $u^{(j)}(t)$ exist everywhere on $\mathbf{R}_0^+$ for $0 \leq j \leq n$ (then guaranteeing the well posedeness of the estimation scheme) so that an extended dynamic 2n-th order system might be built to investigate the closed-loop stability of the adaptive system. The adaptive case for the case when the parameters of (1) are not fully known is investigated in the next section.

## IV. ADAPTIVE CONTROL

### A. Estimation Equations

If the parameters of the system are not fully known then the output trajectory solution may be not calculated via Theorem 1. Thus, equation (8) is replaced with its estimated one which has the same structure but its parameters, the coefficients of the various polynomials in (5) or (8), are replaced by their estimates which are denoted with superscripts "hat" and calculated for all time from an estimation algorithm specified later on. In particular, the output equation (8) is used for filtering purposes what leads to the output estimate:

$$\hat{y}(t) = (F(D) - \hat{M}(D,t))y_f(t) + \hat{\Delta}_0(D,t)u_f(t) + \hat{\Delta}_1(D,t)u_f(t-h) + \hat{\Delta}_2^T(D,t)v_{0f}(t) + \hat{y}_0(t)$$
(26)

The estimation error is calculated directly from the output estimate (26) and the measured output (8)-(10) as:

$$e(t) := y(t) - \hat{y}(t) = -\left(\tilde{\Delta}_0(D,t)u_f(t) + \tilde{\Delta}_1(D,t)u_f(t-h) + \tilde{\Delta}_2^T(D,t)v_{0f}(t) + \tilde{M}(D,t)y_f(t) + \tilde{y}_0(t)\right)$$
(27)

where:

$$\tilde{M}(D,t) := \hat{M}(D,t) - M(D,t) \; ; \; \tilde{\Delta}_i(D,t) := \hat{\Delta}_i(D,t) - \Delta_i(D,t) \, (i=0,2) \; ; \; \tilde{y}_0(t) := \hat{y}_0(t) - y_0(t)$$
(28)

denote parametrical errors. One gets from (26)-(27) that the output may be rewritten by sing estimates as follows:

$$y(t) = (F(D)-M(D,t))y_f(t)+\Delta_0(D,t)u_f(t)+\Delta_1(D,t)u_f(t-h)+\gamma_0(t) \qquad (29.a)$$
$$= (F(D)-\hat{M}(D,t))y_f(t)+\hat{\Delta}_0(D,t)u_f(t)+\hat{\Delta}_1(D,t)u_f(t-h)+e(t) \qquad (29.b)$$

where

$$\gamma_0(t) := \Delta_2^T(D,t)v_{0_f}(t)+y_0(t) \qquad (30)$$

includes the contributions of disturbances including unmodeled dynamics and those of the initial conditions of the filters. The filtered output is the solution to the linear time-varying differential system:

$$\hat{M}(D,t)y_f(t) = \hat{\Delta}_0(D,t)u_f(t)+\hat{\Delta}_1(D,t)u_f(t-h)+e(t) \qquad (31)$$

The adaptive control law maintains the structure of (18)-(20) while replacing the various parameters by their estimates as follows:

$$F_0(D)R(D,h,t)u(t) = T(D)u_c(t)-S(D,h,t)y(t)+L(D)u_0(t) \qquad (32)$$

$$M_0^*(D) = F_0(s)\hat{M}(D,t)R_0(D,t)+\hat{\Delta}_0(D,t)S_0(D,t) \qquad (33.a)$$

$$M_1^*(D,t)-\hat{\Delta}_1(D,t)S_0(D,t) = F_0(D)\hat{M}(D,t)R_1(D,t)+\hat{\Delta}_0(D,t)S_1(D,t) \qquad (33.b)$$

$$L(D)(\hat{\Delta}_0(D,t)u_0(t)+\hat{\Delta}_1(D,t)u_0(t-h))$$
$$+F(D)\hat{M}(D,t)(R_0(D,t)\hat{\Delta}_2^T(D,t)v_0(t)+R_1(D,t)\hat{\Delta}_2^T v_0(t-h))-\hat{\Delta}_1(D,t)\hat{S}_1(D,t)y(t-2h)=0$$
$$(33.c)$$

where $R_i(D,t)$, $S_i(D,t)$ (i=0,1) and $u_0(t)$ are the solutions to the Diophantine equations of polynomials (33.a)-(33.b) below, which exist for all time if the same polynomial degree constraints (22) of Theorem 1 are maintained for all time for their corresponding counterpart estimates. The compensating control input $u_0: \mathbf{R}_0^+ \to \mathbf{R}$ is the solution to the time-varying differential equation (30.c) provided that it exists. The conditions for the existence of such a solution are not very restrictive since it suffices that the estimated functions be everywhere piecewise continuous. Note that, since the estimates are time-varying, the compensating control cannot be now obtained as a convolution integral being the output of a time-invariant transfer function contrarily to the case of known parameters (see (23)-(24)). Eqs. 29, 26, 27 may be rewritten in an equivalent parameterized form which is linear in the parameters as:

$$y(t) = \theta^T \varphi(t)+\gamma_0(t) = \hat{\theta}^T(t)\varphi(t)+e(t) \qquad (34.a)$$
$$\hat{y}(t) = \hat{\theta}^T(t)\varphi(t) \qquad (34.b)$$
$$e(t) = -\tilde{\theta}^T(t)\varphi(t)+\gamma_0(t) \qquad (34.c)$$

where $\hat{\theta}(t)$ is the parameter vector estimate of the parametrical vector $\theta = \left(\theta_y^T, \theta_u^T, \theta_{u'}^T, \theta_{v_0}^T\right)^T$ with $\theta_{v_0} = \left(\theta_{v_{01}}^T, \theta_{v_{02}}^T, ..., \theta_{v_{0n}}^T\right)^T$, which is defined by the coefficients of $(F(D) - M(D))$, $\Delta_i(D)$ (i = 0,1) located as the components of the respective vectors written in the appropriate order so that they are compatible with the regressor vector $\varphi(t)$ is $\varphi(t) = \left(\varphi_y^T(t), \varphi_u^T(t), \varphi_{u'}^T(t)\right)^T$ defined by :

$$\varphi_y(t) = \left(D^{n-1} y_f(t), D^{n-2} y_f(t), ..., y_f(t)\right)^T \qquad (35.a)$$

$$\varphi_u(t) = \left(D^n u_f(t), D^{n-1} u_f(t), ..., u_f(t)\right)^T \qquad (35.b)$$

and $\tilde{\theta}(t) = \hat{\theta}(t) - \theta$, $\forall t \in \mathbf{R}_0^+$ is the parametrical error.

The following assumptions are made on the system (1):

$\theta \in \Omega_i$

**Assumptions: 3**. There exists a known bounded set $\Omega$ which is either a convex compact region or a connected union of a finite number of (disjoint or not) compact sets; i.e. $\Omega = \bigcup_{i=1}^{p} \Omega_i$ such that $\theta \in \Omega_i$; some $i \in \{1, 2, ..., p\}$.

**4**. **(a)** The reference model is stable.

**(b)** The polynomials $M(D)$, $\Delta_i(D)$ (i=0,1) are relatively prime $\forall \theta \in \Omega$, i.e. they have no common zeros when considered as complex functions of D so that the system (1) is a controllable and observable pair.

**(c)** $\int_t^{t+\xi} \left\| \dot{\hat{\theta}}(\tau) \right\| d\tau \leq \kappa_1 \xi + \kappa_2$ for some small $\kappa_1 \in \mathbf{R}_0^+$, some $\kappa_2 \in \mathbf{R}_0^+$ and $\forall \theta \in \Omega$, $\forall t \in \mathbf{R}_0^+$ for any $\xi \in \mathbf{R}^+$.

**5.** $\overline{\delta} := \underset{0 \leq j \leq \text{Max}(n_{\Delta_1}, n_{S_1})}{\text{Max}} \left(\left|\delta_{1j} s_{1j}\right|\right)$ is sufficiently small for all $\theta \in \Omega$ where $\delta_{1j}, s_{1j}$ are the real coefficients of the polynomials $\Delta_1(s)$ and $S_1(s)$. Furthermore

**6**: There exists a non-negative function of time $\gamma : \mathbf{R}_0^+ \to \mathbf{R}_0^+$ satisfying $\gamma_0^2(t) \leq \gamma(t) \leq \varepsilon_1 \underset{0 \leq \tau \leq t}{\text{Sup}} \left(e^{-2\rho_0(t-\tau)} \|\varphi(\tau)\|^2\right) + \varepsilon_2$, $\forall t \in \mathbf{R}_0^+$ for some known real constants $\rho_0 > 0$, $\varepsilon_1 \geq 0$ and $\varepsilon_2 \geq 0$. □

**Remark 3**. Assumptions 3-5 are standard in pole-placement indirect adaptive control algorithms of time-invariant plants and leads to solvability of the diophantine equation associated with the synthesis of the pole-placement based controller in the delay-free nominal case. It means that any delay-free system (1) parameterized in $\Omega$, as well as its associate estimation model, are both controllable and observable. Such a requirement can be easily relaxed and extended to the case when it is stabilizable and detectable since the neglected stable cancellations cause an exponentially decaying neglected term in the control signal what does not modify the properties of the adaptive scheme. The controllability of the estimation model may be guaranteed for all time by using projection of the estimates on the boundary of the $\Omega$-domain, if necessary (see, for instance, [8-10]). Note that Assumption 5 is a version of the third condition of Theorem 3 for the case of known parameters concerning the sufficient smallness of the parameters implying the presence of a double delay h in the closed-loop system to guarantee the stability provided that the two other terms involving zero delay and delay h. respectively, conform a Hurwitz characteristic quasi-polynomial.

**Remark 4**. The projection technique is then incorporated in the estimation algorithm proposed in this paper. However, note that it is not required that neither the zeros of $B(D, e^{-hD})$ nor those associated with the corresponding part of the estimation scheme be stable. Another alternative way which could be used to ensure the controllability of the estimated model for all time is the use of estimates modification procedures when the controllability of the estimated model is lost. The previously proposed modification procedures lead, in general, to high computational costs for delay-free plants of second order or higher (see, for instance [7-8], and references therein) what may lead to implementation difficulties. Assumptions 4-5 are used to guarantee closed-loop stability in the presence of disturbances through the synthesis of a pole-placement based adaptive controller. Finally, Assumption 6 holds if the signal $v_{0f}(t)$ of the disturbance contribution to the state evolution is, in a general problem setting up, the sum of a bounded disturbance and a signal related to the input by a strictly proper exponentially stable function. This is a reasonable assumption in many practical problems and it suffices to know upper-bounds of such constants in order to run the estimation scheme. Assumption 6 may be directly extended to the case when $\gamma(t) \leq \varepsilon_1 \sup_{0 \leq \tau \leq t} \left( e^{-2\rho_0(t-\tau)} \|\varphi(\tau)\|^2 \right) + \varepsilon_2$ still holds with unknown constants $\varepsilon_1$ and $\varepsilon_2$ but estimated by extending the estimation scheme including estimates of such constants as proposed in [7]-[8] for the nonsingular delay-free case. □

The following least-squares type multiple estimation algorithm is proposed. It involves an adaptation relative dead zone which is implemented for closed-loop stabilization under uncertainties:

$$\dot{\hat{\theta}}(t) = P^{1/2}(t) \operatorname{Proj} \{b(t) P^{T/2}(t) \varphi(t) e(t)\}$$
$$\dot{P}(t) = -b(t) P(t) \varphi(t) \varphi^T(t) P(t) \ ; \ P(0) = P^T(0) = P_i^T(0) \geq k_0 I \ (k_0 \in R^+ > 0) \quad (36)$$

for i = 1, 2,..., p ; with Proj {.} being a projection operator, [2], [5-6] used to constraint the estimates of the i-th estimator within the bounded convex region $\Omega_i$ (i=1, 2,....,p) in the light of Assumption 3, and the relative adaptation dead zone being:

$$b(t) = \frac{\alpha_1 s(t)}{1 + \varphi^T(t)P(t)\varphi(t)} \quad \text{some real constant } \alpha_1 > 0 \tag{37.a}$$

$$s(t) = \begin{cases} 0 & \text{if } |e(t)| \leq \vartheta\gamma(t)^{1/2} \\ 1 - \vartheta\frac{\gamma(t)^{1/2}}{|e(t)|}, & \text{otherwise} \end{cases} \tag{37.b}$$

where $\vartheta > 1$ is a design constant.

### B) Properties of the Estimation Algorithm

Note that the estimations of the $\varepsilon$-constants are positive and non decreasing with time until a limit ensured by Theorem 4 below is reached. Such a result is proved in Appendix C, is related to the properties of the estimation algorithm (36)-(37) irrespective of the control law provided that Assumptions 3-6 hold.

**Theorem 4**. The subsequent two items hold:

**(i)** $\|\tilde{\theta}\| \in L_\infty$, $\|\hat{\theta}\| \in L_\infty$

$b(t)(\gamma(t) - e^2(t)) \in L_\infty$, $(b(t)\gamma(t))^{1/2} \in L_\infty$, $b^{1/2}(t)|e(t)| \in L_\infty$

$b(t)(\gamma(t) - e(t))^2 \in L_\infty$ if $\vartheta > 1$ and $\|\dot{\hat{\theta}}\| \in L_2 \cap L_\infty$

**(ii)** $b(t)(\gamma(t) - e(t)) \in L_1 \cap L_\infty$, $b^{1/2}(t)(\gamma^{1/2}(t) + |e(t)|) \in L_1 \cap L_\infty$,
$b^{1/2}(t)|e(t)| \in L_1 \cap L_\infty$
$b^{1/2}(t)(\gamma^{1/2}(t) - |e(t)|) \in L_1 \cap L_\infty$, $(b(t)\gamma(t))^{1/2} \in L_1 \cap L_\infty$, for all and all those signals tend asymptotically to zero as $t \to \infty$. Furthermore,

$\|\dot{\tilde{\theta}}\| \in L_2 \cap L_\infty$ and $\hat{\theta}(t) \to \hat{\theta}_\infty$ $(\|\hat{\theta}_\infty\| < \infty)$ as $t \to \infty$. Also, $b(t)(\gamma(t) - e(t))^2 \in L_\infty$ and
$b(t)(\gamma(t) - e(t))^2 \to 0$ as $t \to \infty$ since $\vartheta > 1$. □

### C) Main Robust Stability Result

The adaptive filtered closed-loop system obtained through the filters of transfer functions $1/F(s)$ is described by the subsequent equations below obtained by filtering (29.b) and (32):

$$\hat{M}(D,h,t)y_f(t) = \hat{\Delta}_0(D,h,t)u_f(t) + \hat{\Delta}_1(D,h,t)u_f(t-h) + e(t) \tag{38.a}$$

$$F_0(D)R(D,h,t)u_f(t)=T(D)u_{cf}(t)-S(D,h,t)y_f(t)+L(D)u_{0f}(t) \quad (38.b)$$

The following technical assumption is made:

**Assumption 7.** $n_{F_0} \geq \text{Max}\left(0, \dfrac{n_T + \ell - n_M}{2}\right).$ □

The above assumption reduces to $n_{F_0} \geq n - n_M$ if $\ell \leq 2n - 2n_M + n_M - n_T = 2(n_1 + n_2) - (n_M + n_T) \leq n_1 + 2n_2 - n_T$. The combined equations (38) may be described through the following auxiliary $2(n_{F_0} + n_M)$-th dynamic system in view of the degree constraints (22.a):

$$\dot{x}(t) = \sum_{j=0}^{1} A_j(t) x(t - jh) + A_2(t) x(t - 2h) + B_0 \bar{u}_{0f}(t) + B_1 \bar{u}_{cf}(t) + e_1 e(t) \quad (39)$$

where $n' = 2n_{F_0} + n_M$ and $e_i$ is the i-th unity Euclidean vector,

$$x(t) = \left( y_f^{(n_M-1)}(t), y_f^{(n_M-2)}(t), \ldots, y_f(t), u_f^{(n'-1)}(t), u_f^{(n'-2)}(t), \ldots, u_f(t) \right)^T \quad (40)$$

The time varying matrix $A_0(t) + A_1(t)$ has constant eigenvalues for all $t \geq 0$ which are the zeros of the quasi-polynomial $M_0^*(s) + M_1^*(s)e^{-hs}$, $\bar{u}_{0f}(t)$ is formed with the compensating signal $u_{0f}(t)$ and its successive time-derivatives up till conform $2n_{F_0} + n_M$ components $A_2(t)$ has in its $(n_M + 1)$-th row the coefficients of the polynomial $\Delta_1(s)R_1(s)$ completed with zeros if necessary until completing a row size of $2(n_{F_0} + n_M)$ while the remaining rows are zero. $B_1$ has all its rows zero except the first one which built with the coefficients of $T(D)$, completed with zeros if necessary. The following closed-loop stability theorem is proved in Appendix C.

**Theorem 5.** The auxiliary extended dynamical system (39)-(40) is globally Lyapunov´s stable if $\varepsilon_1$ and $\bar{\delta}$ are sufficiently small and Assumption 3-6 hold.
As a result, the closed- loop system is globally Lyapunov´s stable. □

## ACKNOWLEDGEMENTS

The author is very grateful to the Spanish Ministry of Education by its partial support of this work through Grant DPI 2006-00714.

## Appendix A. Proofs

**Proof of Theorem 1**. **(i)** The differential equation (2) is obtained directly by combining Eqs. 1. One obtains directly by taking Laplace transforms in (1) and using (6)

$$Y(s) = G(s)U(s) + G'^T(s)(V_0(s) + Ex_0) \quad (A.1)$$

since Assumption 1 implies that M(s) is not identically zero, where $G(s) = G'^T(s)(b + d e^{-hs})$ is the transfer function (8) of the nominal system (1) related to the control input. If Assumption 1 holds, the system (1) is solvable so that the output trajectory $y(t)$ exists for any bounded real initial conditions for all $t \geq 0$. Then, the filtered output and all the remaining filtered signals also exist as solutions of the respective differential system for any given bounded real initial conditions of the corresponding filters. Then, the filtered output trajectory satisfies via (A.1) and (7):

$$Y_f(s) = G_f(s)Y(s) + G_f'^T(s)x_{yF0} \quad (A.2)$$

On the other hand, one obtains for arbitrary initial conditions from the filtered signals (7):

$$U_f(s) = G_f(s)U(s) + G_f'^T(s)x_{uF0} \tag{A.3.a}$$

$$Y_f(s) = G_f(s)Y(s) + G_f'^T(s)x_{yF0} \tag{A.3.b}$$

$$V_{0f}(s) = G_f(s)V_0(s) + G_f'^T(s)x_{v0F0} \tag{A.3.c}$$

$$X_{0f}(s) = G_f'^T(s)(Ex_0) \tag{A.3.d}$$

where identical filters of zero-free transfer functions $G_f(s) = \dfrac{1}{F(s)}$ are introduced with $F(s) = \text{Det}(sI_{n_F} - A_F)$ being a monic polynomial of degree $n_F$, subject to a non-unique factorization:

$$G_f(s) = c_F^T(sI_{n_F} - A_F)^{-1}b_F = G_f'^T(s)b_F \tag{A.4}$$

with $G_f'^T(s) = c_F^T(sI_{n_F} - A_F)^{-1}$, defined by the controllable and observable triple $(c_F^T, A_F, b_F)$ for four identical linear filters of transfer functions $1/F(s)$, where $F(s) = \text{Det}(sI_{n_F} - A_F)$ is a monic polynomial of degree $n_F$ equalizing the order of the real square matrix $A_F$ which generates the filtered signals $u_f(t)$, $y_f(t)$, $v_{0f}(t)$ and $x_{0f}(t)$ with respective initial conditions $x_{yF}(0) = x_{yF0}$, $x_{uF}(0) = x_{uF0}$, $x_{0F}(0) = x_{0F0}$ and $x_{v_0F}(0) = x_{v0F0}$. [Note that $c_F^T \text{Adj}(sI_{n_F} - A_F)b_F = 1$ since the denominator polynomial of the filter is monic and the filter is zero-free then the filter transfer function is cancellation-free (then the state-space realization is controllable and observable]. The substitution of (A.3)-(A.4) into (A.2) and the subsequent use of (A.1) yields the filtered output:

$$\begin{aligned}Y_f(s) &= G_f(s)\Big(G(s)U(s) + G'^T(s)\big(V_0(s) + Ex_0\big)\Big) + G_f'^T(s)x_{yF0} \\ &= G(s)\Big(U_f(s) - G_f'^T(s)x_{uF0}\Big) + G'^T(s)\Big(V_{0f}(s) - G_f'^T(s)x_{v0F0} + G_f(s)Ex_0\Big) + G_f'^T(s)x_{yF0} \\ &= G(s)U_f(s) + G'^T(s)\big[V_{0f}(s) - G_f(s)(Ex_0 + x_{v0F0})\big] + G_f'^T(s)\big[x_{yF0} - G(s)x_{uF0}\big]\end{aligned} \tag{A.5}$$

which is equivalent through (6) to

$$\begin{aligned}&-M(s)Y_f(s) + \big(\Delta_0(s) + \Delta_0(s)e^{-hs}\big)U_f(s) \\ &+ M(s)\Big\{G'^T(s)\big[V_{0f}(s) - G_f(s)(Ex_0 + x_{v0F0})\big] + G_f'^T(s)\big[x_{yF0} - G(s)x_{uF0}\big]\Big\} = 0\end{aligned} \tag{A.6}$$

Note from (A.3b) is equivalent to the linear time-invariant differential equation

$$y(t) = G_f^{-1}(D)y_f(t) = F(D)y_f(t) \tag{A.7}$$

The use of (A.6) into the right-hand side of (A.7) yields directly (8) subject to (8)-(10). □

**Proof of Theorem 3**. **(i)** Note that the closed-loop system (19)-(20) is linear and time-invariant and that the Diophantine equations of polynomials (20.a)-(20.b) are solvable if the unforced system equation (1) is controllable and observable since then its transfer function related to the control input is zero-pole cancellation free. Then, polynomials $R_i(s)$ (monic) and $S_i(s)$ (i=0,1) satisfying the polynomial Diophantine equations (20.a)- (20.b) exist and are unique under the given degree constraints since

(1) $F(s)$ is coprime with $\Delta_i(s)$ for i=0, 1.

(2) The system (1) is controllable and observable which together with the above condition implies that $F(s)M(s)$ is coprime with $\Delta_i(s)$; i=0,1.

(3) $n_F + n_M + n_{R_i} > n_{\Delta_i} + n_{S_i}$ for i=0,1 together with the degree constraints (22.a) imply that the Diophantine equations (22.a)-(22.b) are solvable with polynomials leading to realizable compensators.

(4) The constraint degrees (22.b) for the polynomials L(s) and F(s) imply in view of (22.a) that the control input to the system (1) is sufficiently smooth, so that $u^{(j)}(t)$ exists for all time and all integer $0 \leq j \leq \ell - 1$, and, furthermore, that the compensating signal $u_0(t)$ is obtained from a realizable filter in (20.c).

Therefore, the control injected to the singular uncontrolled system (1), via the control law (18)-(20) subject to the degree constraints (22), is generated with realizable transfer functions and it is sufficiently smooth for its well-posededness. Firstly, assume that the external disturbance and compensating signals $v_0(t)$ and $u_0(t)$ are identically zero for all $t \geq 0$. The closed- loop characteristic equation is:

$$M_0^*(s) + M_1^*(s)e^{-hs} + \Delta_1(s)S_1(s)e^{-2hs} = 0 \tag{A.8}$$

Since $M_0^*(s) + M_1^*(s)e^{-hs} = 0$ has its zeros in $\mathbf{C}_{-\nu}$, it follows from Rouché's theorem for zeros of complex – valued functions that $M_0^*(s) + M_1^*(s)e^{-hs} + \Delta_1(s)S_1(s)e^{-2hs}$ has also its zeros within $\mathbf{C}\setminus\mathbf{C}_{-\nu}^+ \subset \mathbf{C}\setminus\mathbf{C}_0^+$ if for some $\nu_1 \in (0, \nu)$:

$\left|\Delta_1(s)S_1(s)e^{-2hs}\right| < \left|M_0^*(s) + M_1^*(s)e^{-hs}\right|, \forall s \in \gamma_{-\mathbf{C}\setminus\mathbf{C}_{-\nu_1}^+}$ where $\gamma_{-\mathbf{C}\setminus\mathbf{C}_{-\nu_1}^+}$ is a closed simple curve

$\gamma_{-\mathbf{C}\setminus\mathbf{C}_{-\nu_1}^+} := \{s \in \mathbf{C}: \text{Re } s = -\mathbf{C}\setminus\mathbf{C}_{-\nu_1}^+\} \cup \{s \in \mathbf{C}: |s| = \infty \wedge \text{Re } s \leq -\mathbf{C}\setminus\mathbf{C}_{-\nu_1}^+\} \subset \mathbf{C}\setminus\mathbf{C}_0^+$

and, equivalently, if

$$\sup_{\text{Re } s = \nu_1}\left(\left|\frac{\Delta_1(s)S_1(s)}{M_0^*(s)e^{2hs} + M_1^*(s)e^{hs}}\right|\right) < 1 \tag{A.9}$$

As a result, the unforced closed-loop system is globally asymptotically Lyapunov´s stable. If the compensating and reference signals are uniformed bounded for all time then the forced closed-loop system Eqs. 19-20 is globally Lyapunov´s stable as a result.

**(ii)** Note that the compensating signal $u_0(t)$ satisfies the stable time-differential equation:

$$L(D)(\Delta_0(D)u_0(t) + \Delta_1(D)u_0(t-h))$$
$$+ F(D)M(D)(R_0(D)\Delta_2^T v_0(t) + R_1(D)\Delta_2^T v_0(t-h)) - \Delta_1(D)S_1(D)y(t-2h) = \sigma_0(t) \equiv 0$$

, $\forall t \in \mathbf{R}_0^+$ so that the closed-loop equation is that of the reference model. Thus, perfect model-following for all time is achieved under zero initial conditions of (1) and perfect asymptotic model – matching is achieved for any arbitrary set of bounded initial conditions of the uncontrolled system (1) and the reference model, namely, if the compensating signal of controller satisfies (23)-(24).

**(iii)** Note that the transfer function $G_m(s)$ has a cancellation $\Delta_0^+(s)$ so that the reference model is not controllable and observable. However, since the cancellation is stable, the output is identical for all time to that of its minimal state-space realization, i.e. that obtained by removing such a cancellation. $\left(R_0'(s), S_0(s)\right)$ are existing unique polynomials satisfying the Diophantine equation (25.b) whose higher degree right-hand-side polynomials are monic by construction. Since the denominator of the reference model is Hurwitz, its input is uniformly bounded for all time, and the closed-loop system transfer functions equalizes that of the reference model except for the stable cancellation, the perfect model-matching result follows directly. □

### Appendix B. Results About the Weierstrass Canonical Form

**Theorem B. 1**. The following properties hold:

**(i)** If Assumption 1 holds then for any vector function $n + 2n_2$ - tuple $\left(z_0, \Psi([-h, 0], \mathbf{R}), \Psi_{v02}[-h, 0], \mathbf{R}^{n_2}\right) \in A(\Psi)$, the Weierstrass canonical form is controllable at any time $t_1 \geq t + 2h$ and observable at any time $t_1 \geq 0$ independent of the external delay if and only if

$$\text{rank}\left[\alpha_1, W\alpha_1, \ldots, W^{n_1-1}\alpha_1, \beta_1, W\beta_1, \ldots, W^{n_1-1}\beta_1\right] = n_1 \quad (B.1)$$

$$\text{rank}\left[\alpha_2, N\alpha_2, \ldots, N^{n_2-1}\alpha_2, \beta_2, N\beta_2, \ldots, N^{n_2-1}\beta_2\right] = n_2 \quad (B.2)$$

$$\text{rank}\left[\gamma, \overline{W}\gamma, \ldots, \overline{W}^{n-1}\gamma\right] = n \quad (B.3)$$

where $\overline{W} := \text{Diag}(W, N)$ or, equivalently, if and only if

$$\text{rank}\left[\lambda I_{n_1} - W, \alpha_1, \beta_1\right] = n_1 \ ; \ \text{rank}\left[N, \alpha_2, \beta_2\right] = n_2 \quad (B.4)$$

$$\text{rank}\left[\overline{\lambda} I_{n_1} - \overline{W}, \gamma\right] = n = n_1 + n_2 \quad (B.5)$$

$\forall \lambda \in \mathrm{sp}(W)$, the spectrum of W and $\forall \overline{\lambda} \in \mathrm{sp}(W) \cup \{0\}$, respectively.

**(ii)** The transfer function of the nominal system (1) related to the control input is given by

$$G(s) = c^T(sE - A)^{-1}(b + de^{-hs}) = \gamma_1^T(sI_{n_1} - W)^{-1}(\alpha_1 + \beta_1 e^{-hs}) - \sum_{i=0}^{\ell-1} \gamma_2^T N^i (\alpha_2 + \beta_2 e^{-hs}) s^i$$

which is zero-pole cancellation free if and only if the polynomial $M(s) = \det(sE - A)$ has degree $n_M \leq n = n_1 + n_2$ and the Weierstrass canonical form of (1) is controllable and observable; i.e. if and only if $n_M \leq n$ and (B.1)-(B.3) or, equivalently (B.4)-(B.5) hold.

**Proof**. Assume for simplicity and without loss of generality that the system is external disturbance-free. Otherwise, any prefixed final state to be matched under controllability conditions $z^*$ may be correspondingly modified by considering the external disturbance as an extra forcing function so that the controllability conditions are similar to those applicable for the disturbance-free case. Similar considerations apply to observability. Now, note from Popov- Belevitch- Hautus controllability/observability rank tests for the Weierstrass canonical state-space description Eqs. (2) that

$$\mathrm{rank}\left[\alpha_1, W\alpha_1, \ldots, W^{n_1-1}\alpha_1, \beta_1, W\beta_1, \ldots, W^{n_1-1}\beta_1\right] = n_1$$
$$\Leftrightarrow \mathrm{rank}\left[\lambda I_{n_1} - W, \alpha_1, \beta_1\right] = n_1 ; \forall \lambda \in \mathrm{sp}(W) \text{ (equivalently, } \forall \lambda \in \mathbf{C} \text{ since it holds}$$

directly for all complex $\lambda$ which is not an eigenvalue of W),

$$\mathrm{rank}\left[\alpha_2, N\alpha_2, \ldots, N^{n_2-1}\alpha_2, \beta_2, N\beta_2, \ldots, N^{n_2-1}\beta_2\right] = n_2$$
$$\Leftrightarrow \mathrm{rank}\left[\lambda N - I_{n_2}, \alpha_2, \beta_2\right] = n_2 : \forall \lambda \in \mathbf{C}$$
$$\Leftrightarrow \mathrm{rank}\left[N, \alpha_2, \beta_2\right] = n_2 \text{ (since } \mathrm{rank}\left[\lambda N - I_{n_2}, \alpha_2, \beta_2\right] = n_2 \text{ holds for } \lambda = 0 \text{)}$$

and

$$\mathrm{rank}\left[\gamma, \overline{W}\gamma, \ldots, \overline{W}^{n-1}\gamma\right] = n \Leftrightarrow \mathrm{rank}\left[\overline{\lambda} I_{n_1} - \overline{W}, \gamma\right] = n ; \quad \forall \overline{\lambda} \in \mathrm{sp}(\overline{W}) = \mathrm{sp}(W) \cup \{0\}$$

since $\overline{W} := \mathrm{Diag}(W, N)$ and N is nilpotent (then all its eigenvalues are zero). Thus (B.1)-(B.3) and (B.4)–(B.5) are equivalent. For any real $T > T_\varepsilon$ and any fixed $T_\varepsilon \geq 0$, consider the input $u: [T, T+2h) \to \mathbf{R}^2$ defined as $u_T(t) = \begin{cases} \overline{u}_{1T}, & t \in [t, T+h) \\ \overline{u}_{2T}, & t \in [T+h, T+2h] \end{cases}$. Define $\overline{u}_T := (u_{1T}, u_{2T})^T \in \mathbf{R}^2$ and define the sets of admissible initial states $S(\Psi) \subset \mathbf{R}^n$ and admissible initial data $A(\Psi_T) \subset \mathbf{R}^n \times C^{(\ell-1)}([T-h, T], \mathbf{R}) \times C^{(\ell-1)}([T-h, T], \mathbf{R}^{n_2})$, respectively, as follows:

$$S(\Psi_T):=\left\{z(T)=(z_{1T}^T, z_{2T}^T)^T, z_{1T}\in\mathbf{R}^{n_1}, z_{2T}=-\sum_{i=0}^{\ell-1}N^i\left[\alpha_2 u^{(i)}(T)+\beta_2\Psi^{(i)}(T-h)\right]\in\mathbf{R}^{n_2}\right\}$$

$$A(\Psi_T):=\left\{(z_{0T},\Psi):z_{0T}\in S(\Psi_T), \Psi\in C^{(\ell-1)}([T-h,T],\mathbf{R})\right\}$$

for any initialization of a sufficiently smooth control. Then (2.a)-(2.b) become for any bounded $(z_{0T},\Psi)\in A(\Psi_T)$:

$$\dot{z}_1(t)=Wz_1(t)+(\beta_1,\alpha_1)^T(\bar{u}_{2T},\bar{u}_{1T})^T \ ; \ N\dot{z}_2(t)=z_2(t)+(\beta_2,\alpha_2)^T(\bar{u}_{2T},\bar{u}_{1T})^T$$

If $\text{rank}\left[\lambda I_{n_1}-W, \alpha_1, \beta_1\right]=n_1$; $\forall \lambda\in\text{sp}(W)$ and $\text{rank}\left[N, \alpha_2, \beta_2\right]=n_2$ the Weierstrass canonical form is controllable with a piecewise constant control on any time interval $[T, T+2h]$. A similar reasoning would apply to observability. Sufficiency of (B.1)-(B.3) or, equivalently, (B.4)-(B.5) for controllability/ observability has been proved . Necessity might be proved by contradiction. Assume that there is $s_0\in\text{sp}(W)$ such that $\text{rank}\left[s_0 I_{n_1}-W, \alpha_1, \beta_1\right]<n_1$. Then , one gets for $s=s_0$ by taking Laplace transforms in (2.a) for zero initial conditions $\left[s_0 I_{n_1}-W, \alpha_1, \beta_1\right]\left[Z_1^T(s), U(s_0), U(s_0)e^{-hs_0}\right]^T=0$ which holds for some $\left[Z_1^T(s), U(s_0), U(s_0)e^{-hs_0}\right]^T\neq 0$ so that $z_1(t)$ cannot be driven to any given $z(t_1)=z_1^*$ through some control input for any admissible initial data and the system is not controllable. A similar reasoning proofs that that there is no sufficiently smooth control such that the sub-state $z_2(t)$ is controllable for any initial admissible data if $\text{rank}\left[N, \alpha_2, \beta_2\right]<n_2$. In the same way, it is proved that observability fails if $\text{rank}\left[s_0 I_{n_1}-\overline{W}, \gamma\right]<n$ for some $s_0\in\mathbf{C}$. Necessity has been also proved so that Property (i) has been fully proved.

**(ii)** It follows directly from Lemma 1 Eqs. 2.a, 2.c , 3.b and 4.b by taking Laplace transforms the Laplace transform of the output for sufficiently smooth input and admissible initial data is provided that Assumption 1 holds:

$$Y(s)=G(s)U(s)+G_\eta^T(s)\eta(s)+G_I(s) \qquad (B.6)$$

where $G(s)$ is the transfer function related to the control input , $G_\eta^T(s)$ is the transfer matrix related to the disturbance signal vector $\eta(s)=(\eta_1^T(s),\eta_2^T(s))^T$ and $G_I(s)$ is a function being dependent on the initial conditions which are given by:

$$G(s)=c^T(sE-A)^{-1}(b+de^{-hs})=\gamma_1^T(sI_{n_1}-W)^{-1}(\alpha_1+\beta_1 e^{-hs})-\sum_{i=0}^{\ell-1}\gamma_2^T N^i(\alpha_2+\beta_2 e^{-hs})s^i$$

$$(B.7.a)$$

$$G_\eta^T(s) = \left[ \gamma_1^T(sI_{n_1} - W)^{-1}, -\sum_{i=0}^{\ell-1} \gamma_2^T N^i s^i \right] \tag{B.7.b}$$

$$G_I(s) = \gamma_1^T(sI_{n_1} - W)^{-1} z_{10} + \sum_{i=0}^{\ell-1} \sum_{j=0}^{i-1} \gamma_2^T N^i \left(\alpha_2 \Psi^{(j)}(0) + \beta_2 \Psi^{(j)}(-h)\right) s^{i-j-1} \tag{B.7.c}$$

The second identity in (B.6) follows by taking into account (10) and the fact that the transfer function related to the control input is independent of each particular state-space description of (1). The first part of Property (ii) has been proved. Now, note that $n \geq n_M$ ($n = n_M$ holds if and only if the state-space realization is minimal and the system (1) is standard). The proof of the sufficiency part of Property (ii) is direct since if some state-space realization of order n is simultaneously controllable and observable then any other realization of the same order is controllable and observable as well and the transfer function is zero-pole cancellation-free. Now, proceed by contradiction to prove the necessity. Assume that the initial conditions and the external disturbance are zero with no loss in generality so that $Y(s) = G(s)U(s)$. First, note that if $n_M > n$ there are always at least $(n_M - n)$ zero pole cancellations in G(s) so that no state-space realization of G(s) of order $n = n_1 + n_2$ is simultaneously controllable and observable. Now, assume $n_M \leq n$. If there is a single zero-pole cancellation $s = s_0$ of multiplicity k in $G(s)$, of order $n_M$, then one of the following possibilities occur for some integer $k_1$ fulfilling $0 \leq k_1 \leq n_1 \leq k \leq n_M$:

$$\det(sE - A) = s^k \det(sI_{n_1-k_1} - W) \det(sN - I_{n_2+k_1-k}); n_M = k + n_1 - k_1 + \deg(\det(sN - I_{n_2+k_1-k}))$$

if $s_0 = 0$ since N is nilpotent, and

$$\det(sE - A) = (s - s_0)^k \det(sI_{n_1-k} - W) \det(sN - I_{n_2}), \ n_M = n_1 + \deg(\det(sN - I_{n_2}))$$

if $s_0 \neq 0$. Since there is a zero-pole cancellation of multiplicity $k \geq 1$ in G(s) then any controllable and observable realization of G(s) is of order $n - k < n$ what leads to a contradiction. Property (ii) has been proved. □

### APPENDIX C. Proof of Theorem 5

**Proof of Theorem 5**. The unique solution of (39) is

$$x(t) = Z(t,0)x_0 + \int_0^t Z(t,-\tau) \left(A_2(\tau) x(\tau - 2h) + B_0 \bar{u}_{0f}(\tau) + B_1 \bar{u}_{cf}(\tau) + e_1 e(\tau)\right) d\tau \tag{C.1}$$

with $x_0 = x(0)$, $x(\tau) = 0$ for $\tau < 0$, subject to the initial conditions given for the input of (1) and its time-derivatives $u^{(i)} = \Psi^{(i)}(t)$, $\forall t \in [-h, 0]$, $i = 0, 1, \ldots, \ell-1$ and $Z(t,\tau)$ is the evolution operator being the unique solution of

$$\dot{Z}(t,\tau) = \sum_{j=0}^{1} \mathbf{A}_j(t,\tau) Z(t-jh,\tau) \;,\; Z(t,t) = I_{2(n_{F_0}+n_M)},\; Z(t,\tau) = 0,\; \forall \tau > t \in \mathbf{R}$$

(C.2)

From Assumptions 1-6, the unforced linear time-varying dynamic auxiliary system (39)-(40) subject to the evolution operator satisfying (C.2) satisfies $\left\|x^*(t)\right\| = K e^{-\nu_0 t} \left\|x_0^*\right\|$ for some real constants $\nu_0 \in (0,\nu)$ and $K > 0$, since:

a- Assumption 4.c and the smallness of the absolute values of the coefficients of $S_1(D)$ implies also that $\int_t^{t+\xi} \left\|\dot{A}_i(\tau)\right\| d\tau \leq \kappa_{01}\xi + \kappa_{02}$, any $\xi \in \mathbf{R}^+$, some small $\kappa_{01} \in \mathbf{R}_0^+$ (depending on $\kappa_0 \in \mathbf{R}_0^+$), some small $\kappa_{02} \in \mathbf{R}_0^+$, $\forall \theta, \hat{\theta} \in \Omega, \forall t \in \mathbf{R}_0^+, i=0,1$

b- The denominator polynomial of the reference model transfer function is Hurwitz and the time-varying Diophantine Equations associated with the adaptive controller synthesis hold with sufficiently small $\kappa_1$. See [9-10] and Theorem 3.

Since the parametrical error is uniformly bounded for all time from Theorem 4 and, furthermore,

$$|e(t)| = \left|\gamma_0(t) - \tilde{\theta}^T(t)\varphi(t)\right| \leq K_1 \sup_{t-T-2h \leq \tau \leq t}\left(\|x(\tau)\|\right) + \gamma^{1/2}(t) \quad \text{(C.3)}$$

$$\leq K_1 \sup_{t-T-2h \leq \tau \leq t}\left(\|x(\tau)\|\right) + \varepsilon_1^{1/2} \sup_{0 \leq \tau \leq t}\left(e^{-\rho_0(t-\tau)}\|x(\tau)\|^2\right) + \varepsilon_2^{1/2}$$

Now, define disjoint components of $\mathbf{R}_0^+ : \mathbf{I}_{1R} := \{z \in \mathbf{R}_0^+ : |e(t)| \leq \vartheta \gamma(t)^{1/2}\}$ and $\mathbf{I}_{2R} := \mathbf{R}_0^+ \setminus \mathbf{I}_{1R}$. Then,

$$|e(t)| \leq \vartheta \gamma(t)^{1/2}, \forall t \in \mathbf{I}_{1R} \quad \text{(C.4.a)}$$

$$|e(t)| \leq |e(t) - \vartheta \gamma^{1/2}(t)| + \vartheta \gamma^{1/2}(t)$$
$$\leq b^{1/2}(t)|e(t) - \vartheta \gamma^{1/2}(t)| + (1 - b^{1/2}(t))|e(t) - \vartheta \gamma^{1/2}(t)| + \vartheta \gamma^{1/2}(t), \forall t \in \mathbf{I}_{2R}$$

(C.4.b)

$b^{1/2}(t)|e(t) - \vartheta \gamma^{1/2}(t)|$ is integrable on any real interval and $(1 - b^{1/2}(t))|e(t) - \vartheta \gamma^{1/2}(t)|$ is uniformly bounded from Theorem 4, Then,

$$|e(t)| \leq K_1 + \vartheta \gamma^{1/2}(t), \forall t \in \mathbf{R}_0^+ \quad \text{(C.5)}$$

for some $K_1 \in \mathbf{R}_+$. Also,

$$\|A_2(\tau) x(\tau - 2h)\| \leq K_2 \sup_{t-T-2h \leq \tau \leq t}\left(\|x(\tau)\|\right) + K_3 \;;\; \|B_0 \bar{u}_{0f}(t) + B_1 \bar{u}_{cf}(t)\| \leq K_4$$

(C.6)

Then, one gets from (C.1) and (C.3)–(C.4) since the reference control is uniformly bounded for all time

$$\operatorname*{Sup}_{0\leq\tau\leq t}\left(\|x(\tau)\|\right) = K e^{-\nu_0 T}\|x(t-T)\| + K K_2 \operatorname*{Sup}_{t-T-2h\leq\tau\leq t}\left(\|x(\tau)\|\right)\int_{t-T}^{t} e^{-\nu_0(T-\tau)}\,d\tau$$

$$+ K\left(\varepsilon_1^{1/2}\operatorname*{Sup}_{0\leq\tau\leq t}\left(e^{-\rho_0(t-\tau)}\operatorname*{Sup}_{0\leq\tau\leq t}\left(\|x(\tau)\|^2\right)\right) + \varepsilon_2 + K_3 + K_4\right)\int_{t-T}^{t} e^{-\nu_0(T-\tau)}\,d\tau + K K_1$$

(C.7)

Since $\nu_0 > 0$ and $K_2$ and $\varepsilon_1$ are sufficiently small, it follows that $\|x(t)\|$ is uniformly bounded on $\mathbf{R}_0^+$ and then both the auxiliary extended system and the closed-loop one are globally Lyapunov stable. □